\documentclass[12pt]{article}
\usepackage{amsfonts,amsmath,amssymb,bbm}
\usepackage{setspace} 
\def\proof{\noindent{\bf Proof:}\hskip10pt}        
\def\QED{\hfill $\Box$}

\font\tenmath=msbm10 scaled 1200
\font\sevenmath=msbm7 scaled 1200
\font\fivemath=msbm5 scaled 1200
\newfam\mathfam \textfont\mathfam=\tenmath
\scriptfont\mathfam=\sevenmath \scriptscriptfont\mathfam=\fivemath

\begin{document}
\def \\ { \cr }
\def \R{\mathbb{R}}
\def\N{\mathbb{N}}
\def\E{\mathbb{E}}
\def\P{\mathbb{P}}
\def\Z{\mathbb{Z}}
\def\D{\mathbb{D}}
\def\C{\mathbb{C}}
\def\da{^{\downarrow}}
\def \e{{\rm e}}
\def \p{{\cal P}}
\def \s{{\cal S}}
\def \g{{\cal G}}
\newtheorem{theorem}{Theorem}
\newtheorem{definition}{Definition}
\newtheorem{proposition}{Proposition}
\newtheorem{lemma}{Lemma}
\newtheorem{corollary}{Corollary}
\centerline{\LARGE \bf   A second order SDE}
\vskip 2mm
\centerline{\LARGE \bf for the Langevin process reflected}
\vskip 2mm
\centerline{\LARGE \bf    at a completely inelastic  boundary}

\vskip 1cm
\centerline{\Large \bf Jean Bertoin}
\vskip 1cm
\noindent
\centerline{\sl Laboratoire de Probabilit\'es et Mod\`eles Al\'eatoires}
\centerline{\sl Universit\'e Pierre et Marie Curie}
\centerline{\sl and DMA, Ecole Normale Sup\'erieure}
\centerline{\sl Paris, France}
\vskip 15mm

\noindent{\bf Summary. }{\small  It was shown in \cite{Be} that
a Langevin process can be reflected at an energy absorbing boundary.
Here, we establish that the law of this reflecting process can be characterized as the 
unique weak solution to a certain second order stochastic differential equation  with constraints, which is in sharp contrast with a deterministic analog.}
 \vskip 5mm
\noindent
{\bf A.M.S. Classification.}  {\tt Primary 60 H 10 , 60 J 55. Secondary 34 A 12}
\vskip 3mm
\noindent{\bf e-mail.} {\tt jbe@ccr.jussieu.fr}

\begin{section}{Introduction}
Consider the motion of a particle in the half-line $\R_+$ under an external force that governs its acceleration. Assume that the energy of the particle is instantaneously absorbed at  the boundary point $0$, meaning that the velocity of the particle  is always $0$ immediately after hitting $0$. In other words,  the trajectory $(X_t)_{t\geq 0}$ of the particle fulfills the constraints of completely inelastic impacts
\begin{equation}\label{eq3}
\left\{ \begin{matrix}
X_t\geq 0  \,, &\\
X_t=0\ \Rightarrow \ \dot X_{t+} = 0\,.&\\
\end{matrix}\right.
\end{equation} 
and
solves the second order differential equation
\begin{equation}\label{eq1} {\rm d} X_t= \dot X_t {\rm d}t\, , \quad   {\rm d} \dot X_t= F_t {\rm d}t + {\rm d} A_t\quad , \quad A_t =-\sum_{0<s\leq t} \dot X_{s-}{\bf 1}_{\{X_s=0\}}\,,
\end{equation}
where $(F_t)_{t\geq 0}$ denotes the external force.
More precisely, $t\to A_t$ is right-continuous non-decreasing and accounts for the kick induced by the boundary. Specifically, if the particle hits $0$ at time $t$ with incoming velocity $\dot X_{t-}<0$,
then $A_t-A_{t-}=-\dot X_{t-}> 0$ so that $\dot X_{t+}=\dot X_{t-}+(A_t-A_{t-})=0$.

Equation \eqref{eq1} can be viewed as a special case of differential measure inclusions which have been studied initially by Schatzman \cite{Schatzman}; see also Ballard \cite{Ballard} and the references therein.  It is quite remarkable that multiple solutions may exist even in situations when the external force is ${\mathcal C}^{\infty}$.

Following a question raised by Bertrand Maury, we are interested in the case when the external force is a generalized function given by a white noise, i.e. when $F_t=\dot B_t$ with $(B_t)_{t\geq 0}$ a standard Brownian motion.
In this setting, it is natural to consider first the much simpler situation when there is no obstacle at $0$,
that is to introduce  the process with values in $\R$ 
$$Y_t=y_0 + t \dot y_0 + \int_0^t B_s{\rm d}s\,, \qquad t\geq 0.$$
The latter will be called here a {\it free Langevin process}, started from location $y_0\in\R$
and with initial velocity $\dot y_0\in\R$; we refer to Lachal \cite{Lachal} for a rich source of results and references in this area.
It is easily seen that for $y_0=\dot y_0=0$,
$0$ is an accumulation point of the set of times at which the free Langevin process returns to $0$.
Informally, this may suggest that if the energy of the Langevin particle is absorbed at each visit to $0$,
then the particle might never be able to reach a strictly positive velocity, and thus might never take off the boundary. It turns out that this intuition is not correct as we shall see.

 It is convenient to agree that throughout this work, all random processes are implicitly c\`adl\`ag, i.e. their sample paths are right-continuous and possess left-limits everywhere, a.s. In a preceding work \cite{Be}, we established the following result of existence and uniqueness in distribution.

\begin{theorem}\label{Thm1} There exists a strong Markov process $(X_t, \dot X_t)_{t\geq 0}$ with values in $\R_+\times \R$,  started from $X_0=\dot X_0=0$, such that
\begin{equation}\label{eq4}
 {\rm d} X_t= \dot X_t {\rm d}t\ ,\, 
 \int_0^{\infty} {\bf 1}_{\{X_t=0\}}{\rm d}t\,=\, 0\  \hbox{and}\ 
X_t=0\ \Rightarrow \ \dot X_t = 0 
\qquad \hbox{a.s.},
\end{equation}
and which evolves as a free Langevin process as long as $X>0$. Specifically,
for every stopping time $S$ in the natural filtration  of $X$   (after the usual completions),
 if we define
 $\zeta_S=\inf\{t\geq S: X_{t}=0\}$, then
conditionally on $X_S=x_0>0$ and $\dot X_S=\dot x_0$,  the process
$(X_{S+t})_{0\leq t\leq \zeta_S}$
is independent of $ {\mathcal F}_S$ and has the same distribution as
$(Y_t)_{0\leq t\leq \zeta}$, where
$$Y_t=x_0+ t \dot x_0+\int_0^tB_s{\rm d}s \ ,\ \zeta=\inf\{t\geq 0: Y_{t}=0\}$$
and $(B_t)_{t\geq 0}$ is a standard Brownian motion. Further, the preceding requirements determine the distribution of $(X_t, \dot X_t)_{t\geq 0}$.
\end{theorem}
We stress that the strong Markov process $(X_t, \dot X_t)_{t\geq 0}$ has jumps at predictable stopping times (namely, the hitting times by $X$ of the boundary point $0$), and thus fails to be standard; in particular, the Feller property does not hold.

The main purpose of this work is to connect the process characterized in Theorem \ref{Thm1} to Equation \eqref{eq1} when the external force $F$ is a white noise. In this direction, it is convenient to 
rewrite \eqref{eq1} in the form
\begin{equation}\label{eq2} {\rm d} X_t= \dot X_t {\rm d}t\quad , \quad   \dot X_t= B_t + A_t\quad ,\quad A_t =-\sum_{0<s\leq t} \dot X_{s-}{\bf 1}_{\{X_s=0\}}\,.
\end{equation}
We are now able to state the main result of this work.

\begin{theorem}\label{Thm2} {\rm (i)} One can construct on some filtered probability space $(\Omega, ({\mathcal F}_{t})_{t\geq 0}, \P)$ an adapted process
 $(X_t)_{t\geq 0}$  distributed as in Theorem \ref{Thm1}
 and an  $ ({\mathcal F}_{t})$-Brownian motion $(B_t)_{t\geq 0}$,
such that Equations \eqref{eq3}  and \eqref{eq2} hold.

\noindent {\rm (ii)}
Conversely, if on some filtered probability space $(\Omega, ({\mathcal F}_{t})_{t\geq 0}, \P)$,  there is
an $({\mathcal F}_{t})$-Brownian motion  $(B_t)_{t\geq 0}$ and an adapted process $(X_t)_{t\geq 0}$  which satisfies Equations \eqref{eq3}  and \eqref{eq2} and starts with initial conditions $X_0=\dot X_0=0$, then $(X_{t})_{t\geq 0}$ is distributed as in Theorem  \ref{Thm1}.
\end{theorem}

We shall refer to the process $X$ which appears in Theorems \ref{Thm1} and \ref{Thm2} as a {\it Langevin process reflected at a completely inelastic boundary}. Note that we implicitly restrict our attention to the case when the process starts from $0$ with initial velocity $0$, which is obviously the most interesting situation and induces no loss of generality. In some loose sense, Theorems \ref{Thm1} and \ref{Thm2} both state the existence and uniqueness in law of the Langevin process reflected at completely inelastic boundary, but viewed from two different perspectives. Theorem \ref{Thm1} belongs to the framework of the theory of Markov processes and their excursions, whereas Theorem \ref{Thm2} is expressed in terms of  stochastic differential equations.
It is well-known that these two theories  are intimately connected, and one may expect that a soft argument should enable us to deduce 
 Theorem \ref{Thm2} from Theorem \ref{Thm1}. 
 
 In this direction,   the existence of a weak solution to \eqref{eq2} and \eqref{eq3} is rather easy and will be established
in the first part of  Section 2 by investigating, in the framework of stochastic calculus, the explicit construction given in \cite{Be} of the
process specified by Theorem \ref{Thm1}. More precisely, the latter is obtained
from the free Langevin process associated to some standard Brownian motion $(W_t)_{t\geq 0}$ first by a reflection {\it \`a la} Skorohod and then by a non-invertible random time-substitution.

However,  establishing  weak uniqueness in Theorem \ref{Thm2} is less straightforward. 
Indeed, if we aim at applying Theorem \ref{Thm1}, then we have to check {\it a priori} that any weak solution 
$(X_t,\dot X_t)_{t\geq 0}$ to \eqref{eq2} and \eqref{eq3} enjoys the strong Markov property.
But it is well-known that solutions of an SDE fulfill the Markov property only in the situation when weak uniqueness holds for the SDE, and thus Theorem \ref{Thm1} cannot help.
We also stress that weak uniqueness is the most striking aspect of Theorem \ref{Thm2} as it  is  in sharp contrast with the deterministic situation for which (strong) uniqueness can fail even with a smooth forcing.

In the second part of Section 2, we shall observe 
a key point which lies at the heart of the proof of weak uniqueness. From the same Brownian motion $(W_t)_{t\geq 0}$ which is used to construct a weak solution $(X_t,B_t)_{t\geq 0}$, one can also
built another standard Brownian motion $(B'_t)_{t\geq 0}$ which is independent
of $(B_t)_{t\geq 0}$, and such that $(W_t)_{t\geq 0}$
can be recovered from $(X_t,B_t)_{t\geq 0}$ and $(B'_t)_{t\geq 0}$.
Weak uniqueness is established in Section 3. We consider any solution
$(X_t,B_t)_{t\geq 0}$ to \eqref{eq2} and \eqref{eq3} where $(B_t)_{t\geq 0}$ is some Brownian motion. We then introduce an independent standard Brownian motion $(B'_t)_{t\geq 0}$, and using the analysis developed in Section 2, we construct 
from $(X_t,B_t)_{t\geq 0}$ and $(B'_t)_{t\geq 0}$ another Brownian motion $(W_t)_{t\geq 0}$, such that  $(X_t)_{t\geq 0}$ can be recovered from $(W_t)_{t\geq 0}$
in the same way as in Section 2.  

In the final Section, we first make some brief historical comments about the question of uniqueness in the -deterministic- setting of mechanical systems with perfect constraints. For the reader's convenience, we also provide a simple example showing that uniqueness of the solution to \eqref{eq3} and  \eqref{eq1} may fail even when the external force is smooth. Finally, we discuss some open questions regarding strong solutions to \eqref{eq2} and \eqref{eq3}.

\noindent{\bf Nota Bene.} In this paper, I will use the same notation $X,B,W$,... for processes which,  {\it in fine}, will be shown to have  the same distributions.
However the initial definition and assumptions for these processes may be different in different sections. I hope that the reader will find this helpful and not confusing.

\end{section}

\begin{section}{A weak solution}
The first purpose of this section is to check that the construction of Section 2 in \cite{Be} also provides a solution to \eqref{eq2} and \eqref{eq3}. Then we shall study this construction in further details to gain insight for the proof of weak uniqueness.

\subsection{Construction of a weak solution}
We start by recalling the construction of Section 2 in \cite{Be} and some of its properties.

Let $W=(W_t,t\geq 0)$ be a standard Wiener process started from $W_0=0$ and
write $({\mathcal W}_{t})_{t\geq 0}$ for its natural filtration  after the usual completions. Define the free  
Langevin process
$$Y_t\,:=\,\int_0^{t} W_s{\rm d}s\,,\qquad t\geq 0\,,$$
its infimum process
$$I_t:=\inf\{Y_s : 0\leq s \leq t\}\,,$$
and the random  closed set of times when $Y$ coincides with its infimum
$${\mathcal I}:=\{t\geq0: Y_t=I_t\}\,.$$ 

We write ${\mathcal I}^{o}$ for the interior of ${\mathcal I}$ and recall
from Lemma 2 in \cite{Be} that  with probability one, the boundary $\partial {\mathcal I}={\mathcal I}\backslash 
{\mathcal I}^{o}$ has zero Lebesgue measure.
Further the  canonical decomposition of the open set ${\mathcal I}^{o}$ into disjoint open intervals  is given by
$${\mathcal I}^{o}\,=\,\bigcup_{u\in{\mathcal J}}]u,d_u[\,,$$
where
${\mathcal J}$ is the set of times at which $Y$ reaches its infimum
for the first time during some negative excursion of $W$, and $d_u$ the first return time to $0$  for $W$ after the instant $u$. That is
$${\mathcal J}:=\left\{t\geq 0: W_t<0, Y_t=I_t \hbox{ and } Y_{t-\varepsilon}
>I_{t-\varepsilon}\hbox{ for all $\varepsilon>0$ sufficiently small}\right\}$$
and 
$$d_u:=\inf\{s>u: W_s=0\}\,.$$
  It is readily seen that ${\mathcal J}$ can be expressed in the form of a countable family of stopping times in the filtration $({\mathcal W}_{t})_{t\geq 0}$. For instance
${\mathcal J}=\{S_{k,n}: k,n\in\N\}$, 
where $S_{k,n}$ is the $k$-th instant $t$ such that $Y_t=I_t$, 
$Y_{t-\varepsilon}
>I_{t-\varepsilon}$ for all $\varepsilon>0$ sufficiently small
and the velocity at time $t$ fulfills  $\dot Y_t \in ]-1/(n-1), -1/n]$. Note that the $d_{S_{k,n}}$ are then also stopping times.

Finally we introduce the right-continuous time-substitution
$$T_t:=\inf\left\{s\geq 0: \int_0^s {\bf 1}_{\{Y_v>I_v\}}{\rm d}v>t\right\}\,, \qquad t\geq 0\,,$$
and then the free Langevin process reflected at its infimum
 (in the sense of Skorohod) and time-changed by $T_t$, that is we set for every $t\geq 0$
$$X_t:=(Y-I)\circ T_t \,.$$
We mention that the process $t\to Y_t-I_t$ has been studied first by Lapeyre \cite{Lapeyre}.
Clearly,  the process $t\to X_t$  only takes nonnegative values, is continuous, and it can be shown that it possess a right-derivative at every $t\geq 0$ given by
$$\dot X_t= W\circ T_t\,;$$
see Equation (7) in \cite{Be}. We also  set ${\mathcal F}_{t} = {\mathcal W}_{T_{t}}$.

The following proposition establishes the existence stated in Theorem \ref{Thm2}(i).

\begin{proposition}\label{P1} Define 
$$A_t:=
\int_0^{T_t} {\bf 1}_{\{Y_s=I_s\}}{\rm d} W_s\ 
\hbox{ and }\ 
B_t:=\int_0^{T_t} {\bf 1}_{\{Y_s>I_s\}}{\rm d}W_s\,, \qquad t\geq 0\,,$$
so that
$$\dot X_t= A_t+B_t \,, \qquad t\geq 0\,.$$
Then there is the identity
$$A_t =-\sum_{0<s\leq t} \dot X_{s-}{\bf 1}_{\{X_s=0\}}\,, \qquad t\geq 0\,,$$
and $(B_t)_{t\geq 0}$ is an $({\mathcal F}_{t})$-Brownian motion.
As a consequence, $(X_t,\dot X_t)_{t\geq 0}$ is a weak solution to \eqref{eq2} and \eqref{eq3} with initial condition $X_0=\dot X_0=0$.
\end{proposition}

\noindent {\bf Remark : } The fact that the series $\sum_{0<s\leq t} \dot X_{s-}{\bf 1}_{\{X_s=0\}}$ converges for every $t\geq 0$ a.s. can be deduced from Corollary 2 in \cite{Be}. However this fact shall be established directly in the present analysis.

\vskip 2mm

\proof  We express
$\dot X_t=W\circ T_t = A_t+B_t$,
where $A_t$ and $B_t$ are defined in the statement.
The basic facts that have been recalled above imply the identities
\begin{equation}\label{eq5}
\int_0^t {\bf 1}_{\{Y_s=I_s\}}{\rm d} W_s = \int_0^t {\bf 1}_{\{s\in{\mathcal I}\}}{\rm d} W_s =\int_0^t {\bf 1}_{\{s\in{\mathcal I}^{o}\}}{\rm d} W_s
=\sum_{u\in{\mathcal J}}( W_{d_u\wedge t}-W_{u\wedge t})\,.
\end{equation}

On the one hand, the assertion that $(B_t)_{t\geq 0}$ is an $ ({\mathcal F}_{t})$-Brownian motion is seen from the very definition of the time-substitution $T_t$ and the Dambis-Dubins-Schwarz theorem (see e.g. \cite{RY} on its page 181). On the other hand, again by definition, $T_t\not\in {\mathcal I}^{o}$ and
$W_{d_u}=0$ for every $u\in{\mathcal J}$. We deduce from \eqref{eq5} the identity 
$$A_t=-\sum_{u\in{\mathcal J}, u\leq  T_t}W_{u}. $$
Further, it is easily checked that ${\mathcal J}$ coincides with the set of times of the form
$u=T_{s-}$ with $s> 0$ an instant at which $X$ hits the boundary point 
$0$ with a negative incoming velocity (i.e. $X_s=0$ and $X_{s-}<0$). Since $\dot X_{s-}=W\circ T_{s-}$, we conclude that
$$A_t=-\sum_{0<s\leq t} \dot X_{s-}{\bf 1}_{\{X_s=0\}}\,, \qquad t\geq 0\,.$$

To complete the proof, either we observe that if $t$ is an instant at which $X_t=0$, then
$\dot X_{t-}=B_t+A_{t-}$ and thus $$\dot X_t= B_t+A_t=\dot X_{t-}+(A_t-A_{t-})= \dot X_{t-}-\dot X_{t-}=0\,,$$
or we just recall from Equation \eqref{eq4} in Theorem \ref{Thm1} that  $X_{t}=0\Rightarrow \dot X_{t}=0$. \QED

\subsection{Some further properties}

We introduce the  time-substitution
$$T'_t:=\inf\left\{s\geq 0: \int_0^s {\bf 1}_{\{Y_v=I_v\}}{\rm d}v>t\right\}\,, \qquad t\geq 0\,,$$
which can be thought of as the dual to $T_t$.
Next we 
set
$$B'_t:=\int_0^{T'_t}{\bf 1}_{\{Y_v=I_v\}}{\rm d}W_v\,, \qquad t\geq 0\,,$$
and then, for every $x\geq 0$,
$$ \sigma'(x):=\inf\{t\geq 0: B'_t>x\}$$
for the first passage time of $B'$ above the level $x$.

\begin{lemma} \label{L1} With probability one, there is the identity
$$T_t = t + \sigma'(A_t)\,,\qquad t\geq 0.$$
\end{lemma}

\proof 
It will be convenient to write $B'(t):=B'_t$ and observe from the definition of $B'$ 
and $A_t$ (in Proposition \ref{P1}) the identities
\begin{eqnarray*}
\sigma'(A_t ) &=& \inf\left\{\int_0^s{\bf 1}_{\{Y_v=I_v\}}{\rm d}v: B'\left(\int_0^s{\bf 1}_{\{Y_v=I_v\}}{\rm d}v\right) > A_t\right\}\\
&=&  \inf\left\{\int_0^s{\bf 1}_{\{Y_v=I_v\}}{\rm d}v: \int_0^s{\bf 1}_{\{Y_v=I_v\}}{\rm d}W_v > \int_0^{T_t}{\bf 1}_{\{Y_v=I_v\}}{\rm d}W_v\right\}\,.\\
\end{eqnarray*}

Then recall \eqref{eq5}. Observe that for every $u\in{\mathcal J}$, the process
$s\to W_{d_u\wedge s}-W_{u\wedge s}$ is $0$ before time $u$, takes some strictly positive values immediately after time $u$,  reaches its overall maximum for the first time at $d_u$, and remains constant after $d_u$. Note furthermore that the intervals $[u,d_u]$ for $u\in{\mathcal J}$ are pairewise disjoint. It follows that 
whenever $ t\not\in {\mathcal I}^o$, the stochastic integral $s\to \int_0^s{\bf 1}_{\{Y_v=I_v\}}{\rm d}W_v$ attains its overall maximum on the time-interval $[0,t]$
at time $t$, and  if we define $r(t)=\inf\{s>t: s\in{\mathcal I}^{o}\}$, then the first instant when
this stochastic integral exceeds its value at time $t$ is $r(t)$.
Further this stochastic integral remains constant on $[t,r(t)]$.

Applying these observations to the random time $T_t\not\in {\mathcal I}^{o}$,
we conclude that
$$\sigma'(A_t)= \int_0^{r(T_t)}{\bf 1}_{\{Y_v=I_v\}}{\rm d}v =\int_0^{T_t}{\bf 1}_{\{Y_v=I_v\}}{\rm d}v\,,$$
and thus
$$T_t=\int_0^{T_t}{\bf 1}_{\{Y_v>I_v\}}{\rm d}v+\int_0^{T_t}{\bf 1}_{\{Y_v=I_v\}}{\rm d}v
= t + \sigma'(A_t)\,,$$
as it has been stated. \QED

We are now able to establish the following statement, which will provides us with the hint
for establishing weak uniqueness in the next section.

\begin{proposition}\label{P2} The process
$(B'_t)_{t\geq 0}$ is a standard Brownian which is independent
of $(B_t)_{t\geq 0}$. Further, $W$ can be recovered from $(X,B,B')$ as
$$W_{t}= B_{\tau(t)} + B'_{\tau'(t)}\,,$$
where
$$\tau(t):=\inf\{s\geq 0: s+\sigma'(A_{s})> t\}
\ \hbox{ and }\ \tau'(t):= 1-\tau(t).$$
\end{proposition}

\proof That $(B'_t)_{t\geq 0}$ is a Brownian motion which is independent
of $(B_t)_{t\geq 0}$ follows immediately from the definition of $B$ and $B'$ and
Knight's extension of  the Dambis-Dubins-Schwarz theorem (see e.g. \cite{RY} on its page 183). 

Then we simply write
$$W_{t}= \int_0^{t} {\bf 1}_{\{Y_s>I_s\}}{\rm d} W_s + 
\int_0^{t} {\bf 1}_{\{Y_s=I_s\}}{\rm d}W_s=
 B_{\tau(t)}
+ B'_{\tau'(t)}\,,
$$
where 
$$\tau(t):=\int_0^{t} {\bf 1}_{\{Y_s>I_s\}}{\rm d} s
\ \hbox{ and }\ \tau'(t):=\int_0^{t} {\bf 1}_{\{Y_s=I_s\}}{\rm d} s\,.$$
The identity $\tau(t) + \tau'(t)=t$ is obvious. By definition, $T_t=\inf\{s\geq 0: \tau(s)>t\}$
and thus  $\tau(\cdot)$ coincides with the continuous left-inverse of the strictly increasing  time-change $T_{\cdot}$, i.e. $\tau(T_t)\equiv t$. We get from Lemma \ref{L1} that
$\tau(t)=\inf\{s\geq 0: s+\sigma'(A_{s})> t\}$. \QED

\vskip 2mm
\noindent{\bf Remark : }
By a more careful analysis, one could also establish a stronger
result of independence, namely that $X$ and $B'$ are independent processes.
Nonetheless,  as this will not be needed in this work and also follows from the analysis in the next section, we leave the direct proof to the interested reader. In this direction,
we also stress that the Brownian motion $B$ is adapted to the natural filtration of $X$, as one sees from Proposition \ref{P1}. But we do not know whether, conversely, $X$ is adapted to the natural filtration of the Brownian motion $B$, that is whether the solution to \eqref{eq2} is strong. 
\end{section}

\begin{section}{Uniqueness in distribution}
In this Section, we consider some  filtered probability space $(\Omega, ({\mathcal F}_{t})_{t\geq 0}, \P)$. We assume there is
an $({\mathcal F}_{t})$-Brownian motion  $(B_t)_{t\geq 0}$ and an adapted process $(X_t)_{t\geq 0}$  which satisfies Equations \eqref{eq3}  and \eqref{eq2} and starts with initial conditions $X_0=\dot X_0=0$. Our goal is to establish that $(X_t)_{t\geq 0}$
has the distribution of the  process in the preceding Section. In this direction, 
Proposition \ref{P2} points at the role of an independent Brownian motion, so we assume
that the same probability space $\Omega$ can be endowed with another filtration
$({\mathcal F}'_{t})_{t\geq 0}$ such that the terminal sigma-fields 
${\mathcal F}_{\infty}$ and ${\mathcal F}'_{\infty}$ are independent, and that
there exists an $({\mathcal F}'_{t})$-Brownian motion $(B'_t)_{t\geq 0}$.
Clearly, these assumptions induce no loss of generality (as it suffices to enlarge the initial probability space).

Just as in the preceding Section, we then write
$$ \sigma'(x):=\inf\{t\geq 0: B'_t>x\}$$
for the first passage time of $B'$ above level $x\geq 0$, and define 
$$T_t:=t+\sigma'(A_{t})\,,\qquad t\geq 0\,.$$
The process $t\to T_t$ is strictly increasing and thus possesses a continuous left-inverse
$$\tau(t):=\inf\{s\geq 0: s+\sigma'(A_{s})> t\}\,, \qquad t\geq 0\,,$$
i.e. $\tau(T_t)\equiv t$.
Clearly $0\leq \tau(t)\leq t$, and we also set
$$\tau'(t):= t-\tau(t)\,, \qquad t\geq 0\,.$$
Finally we define
$$W_{t}:= B_{\tau(t)} + B'_{\tau'(t)}\,, \qquad t\geq 0\,.$$

The weak uniqueness stated in Theorem \ref{Thm2}(ii) is now a consequence of the following.

\begin{proposition}\label{P3} {\rm (i)} The process $(W_t)_{t\geq 0}$ is a standard Brownian motion.

\noindent{(ii)} The process $(X_{t})_{t\geq 0}$ can be recovered  as
$$X_{t}= (Y-I)\circ T_{t}\,, \qquad t\geq 0\,,$$
where
$$Y_t\,:=\,\int_0^{t} W_s{\rm d}s \hbox{ and } I_t:=\inf\{Y_s : 0\leq s \leq t\}\,.$$

\noindent{(iii)} Finally, there is the identity
$$T_t=\inf\left\{s\geq 0: \int_0^s {\bf 1}_{\{Y_v>I_v\}}{\rm d}v>t\right\}\,, \qquad t\geq 0\,.$$
\end{proposition}

\noindent {\bf Remark : } Proposition \ref{P3} shows that $X$ is distributed as the process which appears in Theorem \ref{Thm1}, and as a consequence, we must have
$\int_0^{\infty}{\bf 1}_{\{X_t=0\}}{\rm d} t=0$ a.s. It may be interesting to point out  that
this property can be checked directly from \eqref{eq2} and \eqref{eq3}. More precisely, 
the set of times $t$ when $X_t=0$ is contained into the zero set of the Brownian semi-martingale $\dot X= B+A$. That the latter has zero Lebesgue measure a.s. can be seen from the occupation density formula for Brownian semi-martingales, see e.g. Corollary 1 in \cite{Protter} on its page 216.

\vskip 2mm

The rest of this Section is devoted to the proof of Proposition \ref{P3};
we start with the first part.

{\noindent{\bf Proof of (i):}\hskip10pt} 
Although one may perhaps establish the result more directly by stochastic calculus,
we shall use an approximation, as this makes the proof more intuitive.
Specifically, pick $\varepsilon>0$
and introduce
$$A_t ^{(\varepsilon)}:=-\sum_{0<s\leq t} \dot X_{s-}{\bf 1}_{\{X_s=0, \dot X_{s-}<-\varepsilon\}}\,, \qquad t\geq 0\,,$$
so that $t\to A_t ^{(\varepsilon)}$ is a non-decreasing process and
$$
\lim_{\varepsilon \downarrow 0}\uparrow A_t ^{(\varepsilon)}
= A_t\,,Ê\qquad t\geq 0\,.
$$
Set also
 $$\tau^{(\varepsilon)}(t):=\inf\{s\geq 0: s+\sigma'(A^{(\varepsilon)}_{s})> t\}\,,$$
 so $t\to \tau^{(\varepsilon)}(t)$ is a continuous non-decreasing process with
 $0\leq \tau^{(\varepsilon)}(t)
\leq t$ and 
$$
\lim_{\varepsilon \downarrow 0}\downarrow \tau^{(\varepsilon)}(t)
= \tau(t)\,,Ê\qquad t\geq 0\,.
$$
Thus, if we define
$$W^{(\varepsilon)}_t:=B_{\tau^{(\varepsilon)}(t)} + B'_{t-\tau^{(\varepsilon)}(t)}\,,Ê\qquad t\geq 0\,,$$
then we have
\begin{equation}\label{eq8}
\lim_{\varepsilon \to 0+} W^{(\varepsilon)}_t= W_t\,.
\end{equation}

Hence, it now suffices  to check that for every $\varepsilon >0$, the process $(W^{(\varepsilon)}_t)_{t\geq 0}$ is a standard Brownian motion. Let us first provide an intuitive explanation. 
 The time-change
$t\to \tau^{(\varepsilon)}(t)$ 
is an absolutely continuous process, and its derivative ${\rm d}{\tau}^{(\varepsilon)}(t)/{\rm d}t:=\dot{\tau}^{(\varepsilon)}(t)$ is a step process which takes alternately the values $1$ and $0$.
The dual time-change $t\to t-\tau^{(\varepsilon)}(t)$ is also absolutely continuous  with derivative $1-\dot{\tau}^{(\varepsilon)}(t)\in\{0,1\}$.
The process  $W^{(\varepsilon)}$ is thus obtained by following alternately the paths of two independent Brownian motions, $B$ and $B'$,
in a way which may remind us of the classical two-arm bandit (see, e.g. \cite{Gittins}). More precisely,
$W^{(\varepsilon)}$ follows $B$ when $\dot{\tau}^{(\varepsilon)}(t)=1$ and follows $B'$ otherwise. The instants when $W^{(\varepsilon)}$ switches from $B$ to $B'$
correspond to the jump times of $A^{(\varepsilon)}$, whereas the instants when $W^{(\varepsilon)}$ switches from $B'$ to $B$ correspond  
to certain first passage times of $B'$. These switching times form an increasing  sequence
of predictable random times for $W^{(\varepsilon)}$, and we can then deduce from the strong Markov property that $W^{(\varepsilon)}$ is a standard Brownian motion.

More precisely, the assumption that $X$ solves \eqref{eq2} and elementary properties
of the free Langevin process easily imply that with probability $1$, the set of times $t$ at which $X_t$ hits the boundary $0$ with incoming velocity $\dot X_{t-}<-\varepsilon$
is both discrete and unbounded. Thus the set of jump times of $A^{(\varepsilon)}$ can be expressed as an increasing sequence of stopping times $J^{(\varepsilon)}_1<J^{(\varepsilon)}_2<\cdots$ where $J^{(\varepsilon)}_0=0$
and
$$J^{(\varepsilon)}_{n}:=\inf\{t>J^{(\varepsilon)}_{n-1}: X_t=0\hbox{ and }
\dot X_{t-}<-\varepsilon\}\,, \qquad n\in\N\,,$$
and  $\lim_{n\to\infty}J^{(\varepsilon)}_n=\infty$.

Write for simplicity $a^{(\varepsilon)}_n=A^{(\varepsilon)}_{J^{(\varepsilon)}_n}$, and
consider the increasing sequence $(\sigma'(a^{(\varepsilon)}_n))_{n\in\N}$.
As each $a^{(\varepsilon)}_n$ is a random variable which is measurable with respect to
${\mathcal F}_{\infty}$ and thus independent of $B'$, the $\sigma'(a^{(\varepsilon)}_n)$
form an increasing sequence of randomized $({\mathcal F}'_t)$-stopping times.
The strong Markov property entails that conditionally on 
$(a^{(\varepsilon)}_n)_{n\in\N}$, the pieces of Brownian paths
$$(B'_{t+\sigma'(a^{(\varepsilon)}_{n-1})}-a^{(\varepsilon)}_{n-1}:{0\leq t < \sigma'(a^{(\varepsilon)}_n)
-\sigma'(a^{(\varepsilon)}_{n-1})})$$
are independent, and for each fixed $n\in\N$, the conditional law of this $n$-th piece is  that of a standard Brownian motion killed when it exceeds $a^{(\varepsilon)}_n-a^{(\varepsilon)}_{n-1}$. 

The process $(W^{(\varepsilon)}_t)_{t\geq 0}$ is obtained by splicing
the sequence of pieces of paths
$$(B_t: 0\leq t <J^{(\varepsilon)}_1)\,,\,
(B'_{t}:0\leq t < \sigma'(a^{(\varepsilon)}_1))\, 
,\, (B_{t+ J^{(\varepsilon)}_1}-B_{ J^{(\varepsilon)}_1}: 0\leq t <J^{(\varepsilon)}_2-J^{(\varepsilon)}_1)\,,\,
\ldots
$$
In particular
$$ W^{(\varepsilon)}_t= \left\{ \begin{matrix}
 B_t  &\hbox{ when } 0\leq t <J^{(\varepsilon)}_1\,,\\
B'_{t-J^{(\varepsilon)}_1}+ B_{J^{(\varepsilon)}_1} &\hbox{ when } 
J^{(\varepsilon)}_1 \leq t < J^{(\varepsilon)}_1 + \sigma'(a^{(\varepsilon)}_1)\,,\\
\end{matrix}\right.
$$
and the strong Markov property of Brownian motion shows that the process
$(W^{(\varepsilon)}_t: 0\leq t < J^{(\varepsilon)}_1 + \sigma'(a^{(\varepsilon)}_1))$
has the same law as $(B_t:0\leq t < \rho_1)$ 
where  $\rho_1$ is the $({\mathcal F}_t)$-stopping time defined as 
$\rho_1:=\inf\{t>J^{(\varepsilon)}_1 : B_t-B_{J^{(\varepsilon)}_1 }>a^{(\varepsilon)}_1\}$.
Splicing more and more pieces, we now see by that an iteration of this argument  that
$(W^{(\varepsilon)}_t: 0\leq t < J^{(\varepsilon)}_n + \sigma'(a^{(\varepsilon)}_n))$
has the same distribution as $(B_t:0\leq t < \rho_n)$ 
where 
$\rho_n:=\inf\{t>J^{(\varepsilon)}_n : B_t-B_{J^{(\varepsilon)}_n }>a^{(\varepsilon)}_n\}$.
Letting $n\to\infty$, we conclude that $(W^{(\varepsilon)}_t)_{t\geq 0}$ is a standard Brownian motion, and an appeal to \eqref{eq8} completes the proof. \QED

Next, let us write 
$${\mathcal D}_A:=\{t>0: A_t>A_{t-}\}$$
 for the
set of times where $A$ is discontinuous.
Observe from \eqref{eq2} that we have also the identification
$${\mathcal D}_A=\{t>0: X_t=0\hbox{ and } \dot X_{t-}<0\}$$
 as the set of instants when
$X$ hits the boundary $0$ with a strictly negative incoming velocity.
An important step in the proof of Proposition \ref{P3} is provided by the following
representation of the processes $\tau(\cdot)$ and $\tau'(\cdot)$.

\begin{lemma}\label{L4} Introduce the random open set
$${\mathcal O}:=\bigcup_{t\in{\mathcal D}_A} ]t+\sigma'(A_{t-}), t+\sigma'(A_{t})[\,,$$
and write ${\mathcal O}^c:=[0,\infty[\backslash{\mathcal O}$ for its complementary set.

\noindent {\rm (i)} The processes $t\to\tau(t)$ and $t\to\tau'(t)$ are both absolutely continuous non-decreasing processes with Stieltjes measures given by
$${\rm d}\tau(t) = {\bf 1}_{{\mathcal O}^c}{\rm d}t\quad , \quad {\rm d}\tau'(t) = {\bf 1}_{\mathcal O}{\rm d}t\,.$$

\noindent {\rm (ii)} We have $W_t\leq 0$ and $X_{\tau(t)}=0$ for every $t\in{\mathcal O}$, a.s.
\end{lemma}

\proof (i) Write $\lambda$ for the Lebesgue measure on $\R_+$. We have to show that 
\begin{equation}\label{eq9}
\lambda\left({\mathcal O}\cap [0,t]\right)= \tau'(t)\,, \qquad t\geq 0\,.
\end{equation}

In this direction, let
${\mathcal A}$ denote the closed range of the process $A_{\cdot}$,
viz. the set of points of the type $A_t$ or $A_{t-}$ for some $t\geq 0$. The complementary set ${\mathcal A}^c:=\R_+\backslash {\mathcal A}$ has a canonical decomposition as union of disjoint open intervals given by 
$${\mathcal A}^c=\bigcup_{t\in{\mathcal D}_A}]A_{t-},A_t[\,.$$
Observe that, since $A_{\cdot}$ is a pure jump process, then for every $t\geq 0$
$$\lambda([0,A_t]\backslash {\mathcal A})=\sum_{s\in{\mathcal D}_A\cap [0,t]}(A_s-A_{s-})= A_t\,,$$
and hence $\lambda({\mathcal A})=0$.

On the other hand, it is well-known that the first passage process $\sigma'$
is a stable subordinator with index $1/2$. In particular, it is purely discontinuous, and, by the L\'evy-It\^o decomposition, the process of its jumps is a Poisson point process.  
Since  ${\mathcal A}$  has zero Lebesgue measure and is  independent of 
 $\sigma'$, ${\mathcal A}$ does not contain any jump time of $\sigma'$, a.s.
 It follows that for every $v\geq 0$
 $$\lambda\left({\mathcal O}\cap [0,v+\sigma'(A_v)]\right)
 = \sum_{s\in{\mathcal D}_A\cap [0,v]}(\sigma'(A_s)-\sigma'(A_{s-}))=\sigma'( A_v)\,.$$
 Recall now that $T_v=v+\sigma'( A_v)$ and observe that 
 $$\tau'(T_v)=T_v-\tau(T_v)=T_v-v=v+\sigma'(A_v)-v=\sigma'(A_v)\,.$$
 We have thus checked that \eqref{eq9} holds for every $t$
 of the form $t=T_{v}$ for some $v\geq 0$, and hence, by approximation, also for 
 every $t$ of the form $t=T_{v-}$ for some $v> 0$.
 
 Finally, suppose that $t\in]T_{v-},T_v[$, where $v=\tau(t)$.
We get from above
$$\lambda\left({\mathcal O}\cap [0,t]\right)
=\lambda\left({\mathcal O}\cap [0,T_{v-}]\right)+ t-T_{v-} = \tau'(T_{v-})+ t-T_{v-}\,.$$
But
$$\tau'(T_{v-})+ t-T_{v-}= T_{v-}-\tau(T_{v-})+t-T_{v-}=t-v=t-\tau(t)=\tau'(t)\,,$$
and we conclude that \eqref{eq9} holds for all $t\geq 0$.

(ii) If $t\in{\mathcal O}$, then $t\in]T_{v-},T_v[$ where $v=\tau(t)$. By definition, we have
$$W_t= B_v+B'_{t-v}\,.$$
On the other hand, $v\in{\mathcal D}_A$ and thus $X_v=0$. Further, by \eqref{eq2} and  \eqref{eq3}, we have $\dot X_v=B_v+A_v=0$.
We deduce that 
$$W_t=B'_{t-v}-A_v\leq 0\,,$$
as $t-v<\sigma'(A_v)$ (because $t<T_v=v+\sigma'(A_v)$). \QED

We are now able to establish  the second part of Proposition \ref{P3}.

{\noindent{\bf Proof of (ii):}\hskip10pt} 
We decompose
$$Y_t=\int_0^t W_s{\rm d}s=\int_0^t W_s{\rm d}\tau(s)
+\int_0^t W_s{\rm d}\tau'(s)\,.$$
The change of variables $s=T_v$ enables us to rewrite the first integral in the sum
as 
$$\int_0^t W_s{\rm d}\tau(s)=\int_0^{\tau(t)} W_{T_v}{\rm d}v
=\int_0^{\tau(t)} (B_{\tau(T_v)}+B'_{T_v-\tau(T_v)}){\rm d}v\,.
$$
Since $\tau(T_v)=v$ and $T_v-\tau(T_v)=T_v-v=\sigma'(A_v)$, the right-hand side
equals 
$$\int_0^{\tau(t)} (B_{v}+A_v){\rm d}v=\int_0^{\tau(t)} \dot X_v{\rm d}v= X_{\tau(t)}\,,$$
where the first equality follows from \eqref{eq2}.

Next, we write $I'_t:=\int_0^t W_s{\rm d}\tau'(s)$,  so that the process
$t\to Y_t-I'_t=X_{\tau(t)}$ is nonnegative.
Further, we know from Lemma \ref{L4}(ii) that $t\to I'_t$ is a non-increasing process
and that the Stieltjes measure ${\rm d}(-I'_t)$ assigns no mass to ${\mathcal O}^c$,
and {\it a fortiori} is supported on the set of times $t$ such that $Y_t-I'_t=X_{\tau(t)}=0$.
An application of Skorohod's reflection principle (see for instance \cite{RY} on its page 239) enables to make the identification $I'_t=\inf\{Y_s: 0\leq s \leq t\} = I_t$. 
We conclude that
$X_t=X_{\tau(T_t)}=Y_{T_t}-I_{T_t}$. \QED

Finally, we turn our attention to the third part of Proposition \ref{P3}.

{\noindent{\bf Proof of (iii):}\hskip10pt} 
On the one hand, we have seen in the proof of part (ii) above that the infimum $I$
of the free Langevin process $Y$ can be expressed as
$$I_t=\int_0^t W_s{\rm d}\tau'(s)=\int_0^t {\bf 1}_{\mathcal O}(s) W_s{\rm d}s\,,$$
where the second identity follows from Lemma \ref{L4}.

On the other hand, we must have $W_t\leq 0$ for every $t$ such that $Y_t=I_t$. Indeed, if we had $W_t>0$ for such a time $t$, then $Y$ would be strictly increasing on some neighborhood of $t$, which is absurd. 
Now for every $t\geq 0$ such that $Y_t=I_t$ and $W_t<0$, $Y$ is strictly decreasing
on some interval $[t,t']$ with $t'>t$, and thus $Y=I$ on $[t,t']$. Since the total time that $W$ spends at $0$ is zero, we also obtain
$$I_t\,=\,\int_0^t {\bf 1}_{\{Y_s=I_s\}}W_s {\rm d}s\,.$$

Using again the fact that the total time that $W$ spends at $0$ is zero, we deduce
by comparison of these two expressions that with probability one, the
random sets ${\mathcal O}$ and $\{s\geq 0:Y_s=I_s\}$ coincide $\lambda$-almost everywhere.  More precisely, recall the notation ${\mathcal I}:=\{s\geq 0:Y_s=I_s\}$ and 
 that the boundary
$\partial{\mathcal I}={\mathcal I}\backslash {\mathcal I}^o$ has zero Lebesgue measure.
We now see that the open sets ${\mathcal O}$ and ${\mathcal I}^o$ coincide a.s.
 As a consequence of Lemma \ref{L4}, 
$$\tau(t)=\int_0^t {\bf 1}_{{\mathcal O}^c}(s) {\rm d}s= \int_0^t {\bf 1}_{\{Y_s>I_s\}} {\rm d}s\,,\qquad t\geq 0\,,$$
and since $T_{\cdot}$ is the right-continuous inverse of $\tau(\cdot)$, this completes the proof. \QED

\end{section}

\begin{section}{Some comments and questions}
We mentioned in the Introduction that  second order differential equations with constraints of the type \eqref{eq3} and \eqref{eq1}  may have multiple solutions even in the situation when the external force $F$ is smooth. This was first pointed out by Bressan \cite{Bressan}, who also made the conjecture that uniqueness holds when the force is a polynomial function of time. 
Schatzman \cite{Schatzman} formulated the general setting of second order differential inclusions, and established a general theorem of existence of solutions.
She also recovered independently Bressan's example of a force of class ${\mathcal C}^{\infty}$ for which such a system possesses multiple solutions.
Percivale \cite{Percivale} was the first to show that uniqueness holds for systems with only one degree of freedom,  when the force is given by an analytic function of the time and  depends neither on the position nor on the velocity of the particle, and then 
Schatzman \cite{Schatzman2} extended this to the much harder case when  the force is an analytic function of time,  position and velocity. 
Finally  Ballard  \cite{Ballard} considered more general discrete systems with several degrees of freedom and established that uniqueness always holds
in the case when the force is analytic. 
We also refer to Maury \cite{Maury, Maury2}, Moreau \cite{Moreau} and Stewart
\cite{Stewart} for numerical schemes for the computation of the motion of bodies systems with inelastic impacts.

For the convenience of the reader, we shall propose here
simple example of an external force of  class ${\mathcal C}^{k}$ (for any fixed $k\in\N$) for which multiple solutions to \eqref{eq3} and \eqref{eq1} exist.
Consider the increasing sequences  $0< \cdots < s_n< t_n< s_{n+1}<\cdots$ given by
$$s_n:=2^{2n}\quad \hbox{and} \quad t_n:=2^{2n+1}\,,\qquad n\in\Z\,.$$
Then introduce the convex increasing  function $\alpha: [0,\infty[\to [0,\infty[$
which is linear on the intervals $[s_n,s_{n+1}]$ and such that $\alpha (s_n)= s_n^{k+3}$.
Similarly, let $\beta : [0,\infty[\to [0,\infty[$
denote   the convex increasing  function which is linear on the intervals $[t_n,t_{n+1}]$ and such that $\beta(t_n)= t_n^{k+3}$. Observe that $\alpha$ and $\beta$ enjoy a property of self-similarity, namely
 $$\alpha (4u)= 4^{k+3}\alpha(u)\quad \hbox{and}\quad \beta (4u)= 4^{k+3}\beta(u)\,,\qquad u>0\,.$$
It should be obvious from a picture that there exists a function $\varphi: ]0,\infty[\to \R$  of class
 ${\mathcal C}^{k+3}$, 
 which is bounded from above by both $\alpha$ and $\beta$, enjoys the same property of self-similarity, viz. $\varphi (4u)= 4^{k+3}\varphi (u)$, and  fulfills the following requirements :
$$
\left\{ \begin{matrix}
& \varphi(u)=\alpha(u)\Leftrightarrow u=s_n \hbox{ for some } n\in\Z,\\
& \varphi(u)=\beta(u)\Leftrightarrow u=t_n \hbox{ for some } n\in\Z,\\
& \dot\varphi(s_n)=\dot\alpha(s_n+)\hbox{ for every }n\in \Z,\\
& \dot\varphi(t_n)=\dot\beta(t_n+)\hbox{ for every }n\in \Z.\\
\end{matrix}\right.
$$ 
 More precisely, one constructs first a function $\varphi$ which satisfies the preceding requirements on the interval $[1,4]$, in such a way that for every $\ell=0,\ldots, k+3$, the $\ell$-th derivative $\varphi^{(\ell)}$
 of $\varphi$ has $\varphi^{(\ell)}(4-)=4^{k+3-\ell}\varphi^{(\ell)}(1+)$. Then $\varphi$ is extended to
 $]0,\infty[$ by self-similarity, and we set $\varphi(0)=0$. Again by self-similarity, we get that
 $\varphi$ is now of class ${\mathcal C}^{k+2}$ on $[0,\infty[$ with $\varphi^{(\ell)}(0+)=0$ for every $\ell=0,\ldots, k+2$.
 The requirements implies that $X_u:=\alpha(u)-\varphi(u)$ solves \eqref{eq3} and \eqref{eq1}
 with $F_u:=-\Ddot\varphi(u)$ and $A_u:=\dot\alpha(u)$. Similarly, $X_u:=\beta(u)-\varphi(u)$ solves \eqref{eq3} and \eqref{eq1}
 with the same external force  and $A_u:=\dot\beta(u)$. Hence Equations \eqref{eq3} and \eqref{eq1}
 have at least two distinct solutions for $F_u:=-\Ddot\varphi(u)$.

Self-similarity  is merely used above as a convenient tool for checking the regularity of the external force at $0$, and a perusal of the argument reveals that  a large class of counter-examples to uniqueness can be built by mimicking the preceding construction, using now an arbitrary strictly convex increasing function $c:[0,\infty[\to [0,\infty[$
($c(u)=u^{k+3}$ in the example above), and arbitrary increasing sequences
$(s_n)_{n\in\Z}$ and $(t_n)_{n\in \Z}$ with no common point and such that
$\lim_{n\to -\infty}s_n=\lim_{n\to -\infty}t_n=0$. The external force $F=-\Ddot \varphi$ may then no longer be smooth; note that in any case, $F$ has strong oscillations near zero, in the sense that $F$ takes negative and positive values at times arbitrarily close to $0$.
This may suggest that, informally, existence of multiple solutions to  \eqref{eq3} and \eqref{eq1} could hold for quite general external forces $F$ with strong oscillations.
In this direction, recall that uniqueness of the solution has only been established for
analytic external forces, see Ballard \cite{Ballard}. 

Theorem \ref{Thm2} is thus in sharp contrast with the preceding observations, even though the uniqueness is only stated there in a weak sense. Hence an important open question is to ask whether {\it pathwise uniqueness} holds for equations \eqref{eq3} and \eqref{eq2}.

Another interesting problem in this vein is to decide whether or not the solution
which has been constructed in Section 2 is adapted to the natural filtration of the Brownian motion $(B_t)_{t\geq 0}$. One says that the solution is {\it strong} in the case when the answer is positive. We refer to Tsirel'son \cite{Tsirelson} for a classical example of an SDE which has a unique weak solution, but no strong solution.

\end{section}

\vskip 1cm
\noindent
{\bf Acknowledgment :} 
I would like to thank Patrick Ballard and Bertrand Maury for useful historical comments and references about the problem which motivated this work.


\begin{thebibliography}{99}

\bibitem{Ballard} P. Ballard : The dynamics of discrete mechanical systems with perfect unilateral constraints {\it Arch. Rational Mech. Anal. \bf 154} (2000), 199-274.

\bibitem{Be} J. Bertoin : Reflecting a Langevin process at an absorbing boundary.
{\it Ann. Probab.} (to appear). Available via {\tt
http://www.imstat.org/aop/future\underline{ }papers.htm}

\bibitem{Bressan} A. Bressan :  Incompatibilit\`a dei teoremi di esistenza e di unicit\`a del moto per un tipo molto comune e regolare di sistemi meccanici, 
{\it Ann. Scuola Norm. Sup. Pisa} Serie III {\bf 14} (1960), 333-348.

\bibitem{Gittins}  J. C. Gittins :
{\em Multi-armed bandit allocation indices.}
Wiley-Interscience Series in Systems and Optimization. John Wiley \& Sons, Chichester, 1989.

\bibitem{Lachal} A. Lachal : 
Applications de la th\'eorie des excursions \`a l'int\'egrale du mouvement brownien.
{\it S\'eminaire de Probabilit\'es  \bf XXXVIII}, Lecture Notes in Math. 1801 (2003), pp. 109-195.

\bibitem{Lapeyre} B. Lapeyre :
Une application de la th\'eorie des excursions \`a une diffusion r\'efl\'echie d\'eg\'en\'er\'ee. {\it Probab. Theory Relat. Fields \bf 87}  (1990), 189-207.

\bibitem{Maury} B. Maury : Direct simulation of aggregation phenomena.
{\it Comm Math. Sci.} supplemental issue No 1  (2004), 1-11.
Available via {\tt http://intlpress.com/CMS/issueS-1/S-1.pdf}

\bibitem{Maury2} B. Maury : A time-stepping scheme for inelastic collisions. 
{\it Numer. Math.} {\bf 102} (2006), 649-679.

\bibitem{Moreau}  J.J. Moreau :
Some numerical methods in multibody dynamics: Application to granular materials. 
{\it Eur. J. Mech., A} {\bf 13}  (1994), 93-114.

\bibitem{Percivale} D. Percivale : Uniqueness in the elastic bounce problem, I,{\it  J. Differential Equations} {\bf 56} (1985), 206-215.

\bibitem{Protter} Ph. Protter : {\em
Stochastic integration and differential equations.} Second edition.
 Applications of Mathematics {\bf 21}. Springer-Verlag, Berlin,  2004. 

\bibitem{RY} D. Revuz and M. Yor : {\em
Continuous martingales and Brownian motion.}
Third edition. Grundlehren der Mathematischen Wissenschaften {\bf  293}. Springer-Verlag, Berlin, 1999.

\bibitem{Schatzman} M. Schatzman : A class of nonlinear differential equations of second order in time. {\it Nonlinear Anal., Theory, Methods Appl.} {\bf 2} (1978), 355-373.

\bibitem{Schatzman2} M. Schatzman : Uniqueness and continuous dependence on data for one dimensional
impact problems.  {\it Math. Comput. Modelling} {\bf 28} (1998), 1-18.

\bibitem{Stewart} D. E. Stewart :
Convergence of a time-stepping scheme for rigid-body dynamics and resolution of Painlev\'e's problem. 
{\it Arch. Ration. Mech. Anal. \bf 145} (1998), 215-260.

\bibitem{Tsirelson} B. S. Tsirel'son :
An example of a stochastic differential equation having no strong solution. {\it Theory Probab. Appl.} {\bf 20} (1975), 416-418; translation from {\it Teor. Veroyatn. Primen.}
{\bf 20} (1975), 427-430. 

\end{thebibliography}
\end{document}